\documentclass[12pt]{article}
\usepackage{amsmath,amssymb,amsthm,latexsym,amscd}

\title{Algebras of distributions \\
of binary semi-isolating formulas \\ for families of isolated
types \\ and for countably categorical theories\footnote{
\newline\indent \ \ \ The work is supported by RFBR (grant
12-01-00460-a).}}
\author{Sergey V.
Sudoplatov\footnote{sudoplat@math.nsc.ru}}
\date{June 05, 2013}
\begin{document}
\maketitle

\begin{abstract}
We apply a general approach for distributions of binary isolating
and semi-isolating formulas to families of isolated types and to
the class of countably categorical theories.

\end{abstract}

\medskip
Algebras and structures associated with isolating and
semi-isolating for\-mu\-l\-as of a theory are introduced in
\cite{ShS, Su122}. We apply this general approach for
distributions of formulas to families of isolated types and to the
class of $\omega$-categorical theories.

\medskip
{\bf Proposition.} {\em For any theory $T$, a nonempty family
$R\subseteq S^1(\varnothing)$ of isolated types, and a regular
family $\nu(R)$ of labelling functions for semi-isolating
formulas, the ${\rm POSTC}_\mathcal{R}$-structure
$\mathfrak{M}_{\nu(R)}$ consists of positive labels and zero, and
each label $u$ has a {\sl complement}\index{Complement!of label}
$\bar{u}$ such that $u\wedge\bar{u}=\varnothing$ and
$u\vee\bar{u}$ is a maximal element. If $R=\{p\}$ then the monoid
$\mathfrak{SI}_{\nu(p)}=\langle M_{\nu(R)},\cdot\rangle$ {\rm
(}for compositions of labels{\rm )} is generated by a Boolean
algebra, for which $u\vee\bar{u}$ correspond to isolating formulas
of $p$.}

\medskip
{\em Proof.} The inclusion $\bigcup\limits_{p,q\in
R}\rho_{\nu(p,q)}\subseteq U^{\geq 0}$ is proved in Proposition
1.1 \cite{Su122}. By the definition, for any pair $(p,q)\in R^2$
of types, each label $v_q$, corresponding to an isolating formula
$\varphi_q(y)$ of type $q(y)$, is a maximal element (among labels
in $\rho_{\nu(p,q)}$). Then for any label $u\in\rho_{\nu(p,q)}$,
the label $\neg u\wedge v_q$ is the complement $\bar{u}$. It
remains to note that for any isolated type $p(x)$, labels in
$\rho_{\nu(p)}$ form a Boolean algebra with the least element
$\varnothing$ and the greatest element $v_p$ such that for any
label $u\in\rho_{\nu(p)}$, $u\wedge\bar{u}=\varnothing$ and
$u\vee\bar{u}=v_p$.~$\Box$

\medskip
By Ryll-Nardzewski theorem and Proposition, we get

\medskip
{\bf Corollary 1.} {\em For any $\omega$-categorical theory $T$, a
nonempty family $R\subseteq S^1(\varnothing)$, and a regular
family $\nu(R)$ of labelling functions for semi-isolating
formulas, the ${\rm POSTC}_\mathcal{R}$-structure
$\mathfrak{M}_{\nu(R)}$ is finite, consists of positive labels and
zero, and each label $u$ has a complement $\bar{u}$.}

\medskip
{\bf Theorem.} {\em For any ${\rm POSTC}_\mathcal{R}$-structure
$\mathfrak{M}$, in which each label is positive or zero and has a
complement, there is a theory $T$, a nonempty family $R\subseteq
S^1(\varnothing)$ of isolated types, and a regular family $\nu(R)$
of labelling functions for semi-isolating formulas such that
$\mathfrak{M}_{\nu(R)}=\mathfrak{M}$.}

\medskip
{\em Proof.} Consider the construction for the proof of Theorem
8.1 \cite{Su122}. We define a family of isolated $1$-types
bijective with the set $\mathcal{R}$ by disjoint unary predicates
${\rm Col}_p$, $p\in \mathcal{R}$. Here we assume that if
$\rho_{\nu(p)}$ consists of one (non\-empty) label, i.~e.,
$\rho_{\nu(p)}=\{0\}$, then $|{\rm Col}_p|=1$, and if
$|\rho_{\nu(p)}|>1$ then ${\rm Col}_p$ contains infinitely many
elements. Besides, we assume that each formula ${\rm Col}_p(x)$
isolates a type marked by $p$. Further scheme is based on a
generic construction coordinated with operations in the ${\rm
POSTC}_\mathcal{R}$-structure $\mathfrak{M}$ and with a
$\unlhd$-ordering (isomorphic to the ordering of labels in
$\mathfrak{M}$) of formulas $\theta_{p,u,q}(x,y)$, witnessing that
realizations of $p$ semi-isolate realizations of $q$,
$u\in\rho_{\nu(p,q)}$, $p,q\in \mathcal{R}$, and forming with
predicates ${\rm Col}_p$ the language of theory under
construction. Since the formulas ${\rm Col}_q(y)$, $q\in
\mathcal{R}$ are isolated, then for any label
$u\in\rho_{\nu(p,q)}$ and any realization $a$ of $p$, we can
define formulas $\theta_{p,u,q}(a,y)$ and
$\theta_{p,\bar{u},q}(a,y)$ so that these formulas complement each
other in ${\rm Col}_q(y)$:
$\theta_{p,u,q}(a,y)\wedge\theta_{p,\bar{u},q}(a,y)\vdash$ and
$\theta_{p,u,q}(a,y)\vee\theta_{p,\bar{u},q}(a,y)\equiv{\rm
Col}_q(y)$. The generic construction allows to get a required
theory with quantifier elimination.~$\Box$

\medskip
Having a finite ${\rm POSTC}_\mathcal{R}$-structure $\mathfrak{M}$
the schema above, generating a theory with quantifier elimination,
produce a required theory being $\omega$-ca\-te\-go\-r\-i\-cal:

\medskip
{\bf Corollary 2.} {\em For any finite ${\rm
POSTC}_\mathcal{R}$-structure $\mathfrak{M}$, in which each label
is positive or zero and has a complement, there is an
$\omega$-categorical theory~$T$, a~nonempty family $R\subseteq
S^1(\varnothing)$, and a regular family $\nu(R)$ of labelling
functions for semi-isolating formulas such that
$\mathfrak{M}_{\nu(R)}=\mathfrak{M}$.}

\end{document}